\documentclass[12pt,aps,onecolumn,nobibnotes]{revtex4}

\usepackage{graphicx}%
\usepackage{dcolumn}
\usepackage{amsmath}
\usepackage{amssymb}

\makeatletter
\def\btt#1{\texttt{\@backslashchar#1}}%
\DeclareRobustCommand\bblash{\btt{\@backslashchar}}%
\makeatother
\begin{document}
\newcommand{\N}{\mathbb{N}}
\newcommand{\R}{\mathbb{R}}
\newcommand{\C}{\mathbb{C}}
\newcommand{\g}{\mathcal{G}}
\newcommand{\A}{\mathcal{A}}
\newcommand{\I}{\mathcal{I}}
\newcommand{\n}{\mathcal{N}_\Theta}

\title[Short Title]{STAR PRODUCTS ON COADJOINT ORBITS}
\author{\underline{M. A. Lled\'o}}
\email{lledo@athena.polito.it}
\affiliation{
 Dipartimento di Fisica, Politecnico di
Torino,\\ Corso Duca degli Abruzzi 24, I-10129 Torino, Italy,
and\\ INFN, Sezione di Torino, Italy.}

\begin{abstract}
We study properties of a family of algebraic star products defined
on coadjoint orbits of semisimple Lie groups. We connect this
description with the point of view of differentiable deformations
and geometric quantization.
\end{abstract}

\maketitle

\section{FAMILY OF DEFORMATIONS OF THE POLYNOMIALS ON THE ORBIT}
\label{family}

Let $G$ be a Lie group of dimension $n$ and $\g$  its Lie algebra.
The Kirillov Poisson structure on the dual space $\g^\star$ is
given by
\begin{equation}
\{f_1,f_2\}(\lambda)=\langle[(df_1)_{\lambda},(df_2)_{\lambda}],\lambda\rangle,
\qquad f_1,f_2 \in C^{\infty}(\g^*), \quad \lambda \in \g^*.
\label{kirillov}
\end{equation} The symplectic leaves of this Poisson structure coincide
with the orbits of the coadjoint action of $G$ in $\g$,
 \[
\langle\rm{Ad}^*(g)\lambda
,Y\rangle=\langle\lambda,\rm{Ad}(g^{-1})Y\rangle \quad \forall\;
g\in G,\quad \lambda\in \g*,\quad Y\in \g. \]

Let $G$ be a compact semisimple group of rank $n$. Then  the
coadjoint orbits are algebraic varieties. Let  $\{p_i\}_{i=1}^m$
be a set of generators of the algebra of $G$-invariant polynomials
on $\g^*$. The coadjoint orbits  are determined by the values of
these polynomials, that is by the equations
\begin{equation}p_i=c_i, \qquad i=1,\dots m.\label{ideal}\end{equation} The regular orbits are those
for which the differentials $dp_i$ are independent \cite{va}. They
are algebraic symplectic manifolds of dimension $n-m$. The ideal
of polynomials vanishing on a regular orbit $\Theta$ is a prime
ideal  generated by the relations (\ref{ideal}), we will denote it
by $\I_0$. The algebra of polynomials on $\Theta$,
\[\rm{Pol}(\Theta)=\rm{Pol}(\g^*)/\I_0\]
is a Poisson algebra.

 A formal deformation of a Poisson algebra $\A$ over $\C$ is a
 $\C[[h]]$-module $\A_h$ which is isomorphic as a module to $\A[[h]]$,
 the  isomorphism $\Psi:\A[[h]]\rightarrow \A_h$ satisfying the
 conditions:

 \noindent {\bf a.} $\psi^{-1}(F_1F_2)=f_1f_2$ mod($h$) where
$F_i\in \A_{[h]}$ are such that $\psi^{-1}(F_i)=f_i$  mod$(h)$,
$f_i\in \A$.

\noindent {\bf b.} $\psi^{-1}(F_1 F_2 - F_2 F_1) = h\{f_1,f_2\}$
mod($h^2$).

In this definition one can substitute $\C[[h]]$ with $\C[h]$. We
will say then that we have a $\C[h]$-deformation. It is clear that
a $\C[h]$-deformation can be extended to a formal deformation,
while the opposite is not in general true.

Given a deformation $\A_h$, one can make the pull back the product
in $\A_h$ to $\A[[h]]$ by the isomorphism $\Psi$. The product
defined in this way is called a star product and is in general
given by a formal series
\[f\star g=\Psi^{-1}(\Psi(f)\Psi(g))=fg +\sum_{n>1}h^nB_n(f,g) \]
where $B_n$ are bilinear operators. If $\A$ is some space of
functions and $B_n$ are bidifferential operators we say that the
star product is differential. It follows that the star product can
be extended to the whole space of  $C^\infty$ functions, but only
as a formal deformation \cite{rr}.  By choosing another
isomorphism $\Psi'$, one could obtain a star product that is not
differential. So a star product that is not differential can be
isomorphic to a star product that is differential. We will see
examples of this situation later.

A formal (and $\C[h]$) deformation of $\rm{Pol}(\g^*)$ is given by
the enveloping algebra $U_h$ of the Lie algebra  with the bracket
 $h[\,\cdot \,,\,\cdot\,]$, where $[\,\cdot \,,\,\cdot\,]$ is the bracket on $\g$.
  (The tensor algebra needs
to be taken over $\C[[h]]$). One choice for $\Psi$ is the  Weyl
map or  symmetrizer. If $x_1,\dots, x_n$ are coordinates on $\g^*$
and $X_1,\dots ,X_n$ are the corresponding generators of $U_h$,
the Weyl map is
\[W(x_{i_1}\cdots x_{i_p})=\frac{1}{p!}\sum_{s\in
S_p}X_{i_{s(1)}}\cdots X_{i_{s(p)}}.\] The star product
\[f\star_S g=W^{-1}(W(f)W(g))\] can be expressed in terms of
bidifferential operators, so it can be extended to the whole
$C^\infty(\g^*)$.

$\star_S$ is not tangential to the orbits, so it cannot be
restricted to one of them. Nevertheless, the formal deformation
$U_{[h]}$ can be used to induce a deformation of
$\rm{Pol}(\Theta)$. This was developed in \cite{fl}. The idea is
to find an ideal $\I_h$ such that the diagram
\begin{eqnarray*} \rm{Pol}(\g^*)&\longrightarrow&U_{[h]}\\
\Big\downarrow& &\Big\downarrow\\
\rm{Pol}(\Theta)&\longrightarrow&U_{[h]}/\I_{[h]}
\end{eqnarray*}
commutes. The vertical arrows are the natural projections, the
horizontal ones indicate deformations. The ideal $\I_h$ is
generated by
\[W(p_i)-c_i(h)=P_i-c_i(h),\quad c_i(0)=c_i^0,\qquad i=1,\dots n.\]
The ideal is $\rm{Ad}_G$-invariant since $P_i$ are Casimirs of
$U_{[h]}$, so there is a natural action of $G$ on
$U_{[h]}/\I_{[h]}$.  The same construction works with $\C[h]$. We
will consider only $c^i(h)$ such that its degree in $h$ is not
bigger than the degree of $p_i$. In this context one can show that
$\I_h$ is a prime ideal \cite{fll}. Also, the algebras can be
specialized to a value of $h$, say $h_0$, by quotienting with the
proper ideal generated by $h-h_0$. Analyzing the representations
of the specialized algebras, one can see that in general, they are
not isomorphic for ideals with $c_i(h)\neq c'_i(h)$.

\section{STAR PRODUCTS ON THE POLYNOMIALS ON THE ORBIT}

We consider the example of $S^2$ for clarity, although the
argument can be extended to all other compact orbits \cite{fl,
fll}. $\g=\rm{su}(2)$ has dimension 3 with basis $\{X,Y,Z\}$, \[
[X,Y]=Z,\quad [Y,Z]=X,\quad [Z,X]=Y.\] The unique invariant
polynomial is $p(x,y,z)=x^2+y^2+z^2$ and  the Casimir element is
$P=X^2+Y^2+Z^2$. The regular orbits are given by
\[x^2+y^2+z^2=c\]
for $c>0$. A basis of $\rm{Pol}(S^2)$ is $\{[x^my^nz^\nu],\;
m,n\in \N,\; \nu=0,1\}$.  An isomorphism $\rm{Pol}(S^2)[h]\approx
U_h/\I_h$, is given by
\[\tilde\Psi([x^my^nz^\nu])= [X^mY^nZ^\nu],\]
since $\{[X^mY^nZ^\nu],\; m,n\in \N,\; \nu=0,1\}$ is a basis of
$U_h/\I_h$. Define the isomorphism
$\Psi:\rm{Pol}(\rm{su}(2)^*)[h]\rightarrow U_h$
\begin{eqnarray*}&\Psi(x^my^nz^\nu)=X^mY^nZ^\nu, \qquad &m,n\in \N\\
&\Psi(x^my^nz^r(p-c^0))=X^mY^nZ^r(P-c(h)), \qquad &m,n\in \N
;\end{eqnarray*} which sends the ideal $\I_0$ into the ideal
$\I_h$, so it passes to the quotient, where it gives the
isomorphism $\tilde \Psi$. The corresponding star product on
$\rm{Pol}(\rm{su}(2)^*)[h]$ restricts to $\rm{Pol}(S^2)$.

This star product is not differential, as it is shown in
\cite{fll}, but it is isomorphic to $\star_S$. In addition, for an
orbit in a neighborhood of this one, $p-c^0-\Delta c^0=0$, $\Psi$
doesn't preserve the ideal.

Another way of giving a basis is using the decomposition
\[\rm{Pol}(\g^*)\approx I\otimes H,\]
where $I$ is the algebra of invariant polynomials and $H$ is the
space of harmonic polynomials, $H\approx \rm{Pol}(\Theta)$. We
define the  isomorphism $\Phi:\rm{Pol}(\rm{su}(2)^*)[h]\rightarrow
U_h$
\[\Phi((p-c)^m\otimes \eta_m)=(P-c(h))^m\tilde\Phi(\eta_m),\qquad \eta_m\in H\]
where $\bar\Phi$ is any isomorphism
$\tilde\Phi:\rm{Pol}(S^2)[h]\rightarrow U_h/\I_h$. A star product
of this kind was first written down in Ref. \cite{cg}, where
$\tilde\Phi$ was chosen in terms of the Weyl map,
\[\tilde\Phi([\eta])=[W(\eta)],\]
and $c(h)=c$. We will denote this product by $\star_P$. It has the
nice properties that it restricts to all the orbits in a
neighborhood of the regular orbit and that it is ``covariant",
\[gf_1\star gf_2=g(f_1\star f_2).\]
Nevertheless, it is not differential, as it was shown in
\cite{cg}.

Finally, it was proven in \cite{fl} that $U_\hbar/\I_\hbar$, with
$c(\hbar)=l(l+\hbar)$ corresponds to the algebra of geometric
quantization  in the formalism of Ref.  \cite{vo}

\section{DIFFERENTIAL AND TANGENTIAL STAR PRODUCTS}

In this section we want to consider differential star products on
$\g^*$ and on $\Theta$, and to see the relation with the algebraic
approach of the previous section. In Ref. \cite{ko}, the
differential deformations of a Poisson manifold $X$ modulo gauge
equivalence are shown to be in one to one correspondence with the
formal Poisson structures
\[\alpha=h\alpha_1 +h^2\alpha_2 +\cdots, \qquad [\alpha,\alpha]=0,\]
( $\alpha_i$ are bivector fields and $[\,\cdot\,,\,\cdot\,]$
 is the
Schouten-Nijenhuis bracket) modulo the action of formal paths in
the diffeomorphism group. So for every Poisson structure $\beta$,
one can associate canonically an equivalence class of star
products, the one corresponding to $h\beta$. If there are formal
structures starting with $h\beta$ which are not equivalent to
$h\beta$ through a diffeomorphism path, then one has star products
not equivalent to the canonical one such that
\[f\star g-g\star f=h\beta(f,g)\;\;\rm{mod}(h).\]

In the case of symplectic manifolds, these structures are
classified by $H^2(X)[[h]]$. Since the compact coadjoint orbits
have non trivial second cohomology group, we have more than one
equivalence class of differential star products with term of first
order the same Poisson bracket.

In the case of $\g^*$ with the Kirillov Poisson structure, it
depends on the Lie algebra cohomology of $\g$. So for a semisimple
Lie algebra there is only one equivalence class \cite{fll}.
$\star_S$ is a representative of this equivalence class. It is not
tangential to the orbits, and in fact, it was shown in Ref.
\cite{cgr} that no tangential star product could be extended over
0 for a semisimple Lie algebra.

Nevertheless, a regular orbit has always a neighborhood that is
regularly foliated, $\n\approx\Theta\times \R^m$. Since the
Poisson structure is tangential, the coordinates on $\R^m$ can be
considered as parameters, so one has in fact a family of Poisson
structures on $\Theta$ smoothly varying with the parameters $p_i$,
$\beta_{p_1,\dots p_m}$. Kontsevich's construction of the
canonical star product  gives a star product smoothly varying with
the parameters $p_i$, or, interpreting it in the other way, a
tangential star product canonically associated to $\beta$. It
follows that $\star_S$, when restricted to $\n$, is equivalent to
a tangential star product. We denote it by $\star_T$.

We have three different products,

\smallskip

\noindent  $\star_S$. It is differential, not tangential and
defined on $\g^*$.

\noindent  $\star_P$. It is not differential,  tangential and
defined on $\g^*$. (The ideal has been chosen so that
$c_i(h)=c_i^0$).

\noindent  $\star_T$. It is differential, tangential and defined
only on $\n$.

\smallskip

$\star_S$ restricted to the polynomials  is isomorphic to
$\star_P$. There is then an algebra  homomorphism
\[\varphi:(\mbox{Pol}(\g^*)[[h]],\star_P)\rightarrow(C^\infty(\g^*)[[h]],*_S).\]
We have that \[\varphi(p_i-c_i^0)=p_i-c_i^0,\] so $\I_0\subset
\mbox{Pol}(\g^*)[[h]]$ is sent by $\varphi$ into $\I_0\subset
C^\infty(\g^*)[[h]]$.

Restricting $\star_S$ to $\n$ we have an algebra homomorphism
\[\rho:(C^\infty(\n)[[h]],\star_S)\rightarrow(C^\infty(\N)[[h]],*_T).\]
It is not difficult to see that although in general $\rho$ doesn't
send polynomials into polynomials, the algebra homomorphism
structure implies that  \cite{fll} \[\rho(p_i-c_i^0)=p_i-c_i^0.\]

By composing $\rho\circ\varphi$, one obtains an homomorphism from
a non differential star product to a differential one, such that
both star products are tangential and the ideal $\I_0$ is mapped
into the ideal $\I_0$. The homomorphism passes to the quotient, so
the algebraic star product described in Section 2 is shown to be
homomorphic to the differentiable star product associated by
Kontsevich's map.

We note that we have chosen an algebraic star product with
$c_i(h)=c_i$. This star product is not the one obtained from
geometric quantization. The differential approach to quantization
and geometric quantization, although they have similar features in
the case of $\R^{2n}$ \cite{gv}, seem not to give for compact
coadjoint orbits, homomorphic algebras.

\
\end{document}